\documentclass[a4paper, 14pt]{article}
\usepackage{amsfonts,amssymb,amsmath}
\usepackage{hyperref} 
\usepackage[T2A]{fontenc} 
\usepackage[cp1251]{inputenc} 
\usepackage[russian]{babel}

\def\qed{ \ \hfil$\square$}

\def\Ker{{\hskip0.3mm\rm  Ker\hskip0.5mm}}
\def\Im{{\hskip0.3mm\rm  Im\hskip0.5mm}}

\makeindex

\newcommand{\Hom}{{\rm Hom}}

\newcommand{\ord}{{\rm ord}}
\newcommand{\SL}{{\rm SL}}

\newcommand{\lcm}{{\rm LCM~\!\!}}

\newcommand{\Q}{{\mathbb Q}}

\newcommand{\Z}{{\mathbb Z}}

\newcommand{\Qb}{\overline{\mathbb Q}}

\newcommand{\Ar}{{\mathcal A}}

\newcommand{\Br}{{\mathcal B}}

\newcommand{\Lr}{{\mathcal L}}

\newcommand{\Sr}{{\mathcal S}}

\newcommand{\al}{\alpha}
\newcommand{\De}{\Delta}

\newcommand{\Ga}{\Gamma}

\newcommand{\La}{\Lambda}

\font\teneusm=eusm10 \font\seveneusm=eusm7 
\font\fiveeusm=eusm5 
\newfam\eusmfam 
\textfont\eusmfam=\teneusm 
\scriptfont\eusmfam=\seveneusm 
\scriptscriptfont\eusmfam=\fiveeusm

\def\mat #1,#2,#3,#4,{\left({#1\atop #3}{#2\atop #4}\right)}
\def\bra#1,{{\left\lbrace {#1}\right\rbrace}}

\def\si{\sigma}

\def\l1{\langle}

\newcommand{\B}{\left(\begin{array}{cc}}
\newcommand{\E}{\end{array}\right)}

\def \fns{{${}^{*)}$}}

\newcommand{\comm}[1]
{\fns\marginpar{$\boxed
{\hskip-6pt
{\small {\sf 
\begin{tabular} {l}
 #1
\end{tabular}
}
}
}
$
}
}
\def \?  {\comm{check ?}}

\usepackage{euscript}
\let\scr=\EuScript

\let\mathcal=\scr           
\def\ang#1,{{\left\langle {#1}\right\rangle}} 
\newcommand{\ds}{\displaystyle}

\def\qed{ \ \hfil$\square$}

\def\Ker{{\hskip0.3mm\rm  Ker\hskip0.5mm}}
\def\Im{{\hskip0.3mm\rm  Im\hskip0.5mm}}

\font\teneusm=eusm10 \font\seveneusm=eusm7 
\font\fiveeusm=eusm5 
\newfam\eusmfam 
\textfont\eusmfam=\teneusm 
\scriptfont\eusmfam=\seveneusm 
\scriptscriptfont\eusmfam=\fiveeusm

\font\tengothic=eufm10
\font\sevengothic=eufm7
\font\fivegothic=eufm5
\newfam\Gothic
\textfont\Gothic=\tengothic
\scriptfont\Gothic=\sevengothic
\scriptscriptfont\Gothic=\fivegothic

\def\mat #1,#2,#3,#4,{\left({#1\atop #3}{#2\atop #4}\right)}
\def\bra#1,{{\left\lbrace {#1}\right\rbrace}}

\def\si{\sigma}

\def\l1{\langle}

\def\Numero{${\rm N\sp\circ}$}

\let\scr=\EuScript

\let\mathcal=\scr           

\let\bold=\boldsymbol

\def\bpsi{{\bold \psi}}

\def\vin{{ {\tiny \mid }  
\kern-7.29pt 
\bigcup }}

\def\ang#1,{{\left\langle {#1}\right\rangle}} 

\def \cds{{\cdot{\dots}\cdot}}

\normalsize

\newcommand{\HH}{\mathbb H}

\newcounter{ncours}{\setcounter{ncours} {1}}

\def\SL {\mathop{\rm SL}\nolimits}

\begin{document}
\title{Тройные произведения семейств  Колмана}

\author{А. А. Панчишкин \\
 {\it Памяти дорогого Кости Бейдара
}}

\date
{ 
  Institut Fourier, Universit\'e Grenoble-1\\
  B.P.74, 38402 St.--Martin d'H\`eres, FRANCE\\
  e-mail : panchish$@$mozart.ujf-grenoble.fr, FAX:  33 (0) 4 76 51 44 78
 {\tt http://www-fourier.ujf-grenoble.fr/\~{}panchish/30LVM}
}

\maketitle
\tableofcontents
\begin{abstract}
{ 
Модулярные формы изучаются с точки зрения компютерной алгебры, 
а также как элементы $p$-адических Банаховых модулей.

Представлены методы решения проблем теории чисел посредством производящих функций и их связи с модулярными формами. 

В частности, обсуждаются специальные значения $L$-функций.

Для простого,числа $p$ рассматриваются тройки классических модулярных форм
\begin{align*}
 f_j(z)=\sum_{n=1}^\infty a_{n,j}e(nz)\in \Sr_{k_j}(N_j, \psi_j),\  (j=1, 2,3)
\end{align*}
весов  $k_1, k_2, k_3$, уровней  $N_1, N_2, N_3$,и характеров  $\psi_j \bmod N_j$.

Описаны $p$-адические $L$ функции четырёх переменных, связанные 
с тройными произведениями  семейств  Колмана
$$
k_j{}\mapsto \left \{f_{j,k_j{}}= \sum_{n=1}^\infty a_{n,j}(k{}) q^n\right \} 
$$  
параболических форм положительного наклона
$\sigma_j=v_p(\alpha_{p, j}^{(1)}(k_j{}))\ge 0$
где $\alpha_{p,j}^{(1)} = \al_{p,j}^{(1)}(k_j{})$
собственные значения оператора Аткина $U=U_p$.
\footnote{
Ключевые слова: модулярные формы, 
$p$-адические $L$-функции,
$p$-адические банаховы модули.

УДК: Теория чисел \rm (MSC 11F60)}
}
\end{abstract}

\section{\bf Введение}

{\bf Модулярные формы как объекты компьютерной алгебры }

Костя Бейдар был очень жизнерадостным
человеком, про которых говорят : "душа компании". 
Oн любил танцевать, но также и 
любил решать открытые математические проблемы, и 
он уважал комьютерную алгебру.
В некоторых работах он использовал 
$p$-адические числа, Банаховы кольца и модули.

Решение проблемы Колмана-Мазура о существовании $p$-адических 
$L$-функций двух переменных, связанных с собственными семействами 
положительного наклона было дано автором в предыдущей статье 
\cite{PaTV} (Invent. Math. v. 154, N3 (2003), pp. 551 - 615).

Цель данной работы -- перенести результаты статьи \cite{PaTV} 
на тройные пройзведения Гарретта семейств  модулярных форм Колмана.

\medskip
\begin{tabular}{ll}
\begin{tabular}{l}  
Мы рассматриваем \\
модулярные формы как \\
1) {\it степенные ряды} \\ 
$\ds f=\sum_{n=0}^\infty a_n q^n\in {\mathbb C}[[ q]]$ и как \\
\hskip-0.25cm\begin{tabular}{l}
2) {\it голоморфные функции} \\  
{\it на верхней полуплоскости} \\
$ \HH = \{ z\in {\mathbb C}\ |\ \Im\ z>0\}$\\
\end{tabular}
\end{tabular}
\hskip0.3cm
\begin{tabular}{|l} 
где $q=\exp(2\pi i z)$, \\
$z\in \HH$,  и рассмотрим\\
$L$-функцию \\ 
$\ds L(f, s,\chi)=\sum_{n=1}^\infty \chi(n)a_n n^{-s}$ \\
\it для любого характера Дирихле\\
$\chi:(\Z/N\Z)^*\to {\mathbb C}^*$.
\end{tabular}
\end{tabular}

\medskip
\noindent
{\bf Знаменитый пример - функция Рамануджана $\tau (n)$
}

\begin{tabular}{ll} 
\begin{tabular}{l|}  
{\it Функция $\Delta$ (переменной $z$)}\\
{\it определена формальным выражением} \\
$\Delta=\sum_{n=1}^\infty  \tau (n) q^n$\\
\\
$=q\prod_{m=1}^\infty (1-q^m)^{24}= q-24q^2+252q^3+\cdots$\\
(\it модулярная форма\\
{\it относительно группы $\Ga = \SL_2(\Z)$}).
\end{tabular}
\begin{tabular}{l} 
$\tau (1)=1,  \tau (2)=-24$,\\
$\tau (3)=252, \tau(4)=-1472$\\
$\tau (m) \tau (n) =\tau (mn)$ \\
для $(n,m)=1$, \\
$|\tau(p)|\le 2p^{11/2}$\\
{\sf(Ramanujan-Deligne)}\\
{\it для всех простых чисел}
$p$ .
\end{tabular}
\end{tabular}

\medskip
\noindent
{\bf Быстрое вычисление функции Рамануджана:}

Положим
$ \ds
h_k:=\sum_{n=1}^\infty \sum_{d|n}d^{k-1}q^n
=\sum_{d=1}^\infty \frac {d^{k-1}q^d}{1-q^d}.
$ 

Доказывается: 
$\De=(E_4^3-E_6^2)/1728$ где
$E_4= 1+240h_4$ и
$ E_6= 1 - 504h_6$: 

\medskip
\noindent
{\bf Вычисление с PARI-GP}

\noindent
(см. \cite{BBBCO}, The PARI/GP number theory system. 
{\tt http://pari.math.u-bordeaux.fr)}

$ \ds
h_k:=\sum_{n=1}^\infty \sum_{d|n}d^{k-1}q^n
=\sum_{d=1}^\infty \frac {d^{k-1}q^d}{1-q^d}
\Longrightarrow
$ 

\begin{verbatim}
gp > h6=sum(d=1,20,d^5*q^d/(1-q^d)+O(q^20))
gp > h4=sum(d=1,20,d^3*q^d/(1-q^d)+O(q^20)
gp > Delta=((1+240*h4)^3-(1-504*h6)^2)/1728
\end{verbatim}

\begin{verbatim}
q - 24*q^2 + 252*q^3 - 1472*q^4 + 4830*q^5 - 6048*q^6 - 16744*q^7 
+ 84480*q^8 - 113643*q^9 - 115920*q^10 + 534612*q^11 
- 370944*q^12 - 577738*q^13 + 401856*q^14 + 1217160*q^15 
+ 987136*q^16 - 6905934*q^17+ 2727432*q^18 + 10661420*q^19 + O(q^20)
\end{verbatim}

\noindent
{\bf Сравнение Рамануджана: }
$\ds \tau (n) \equiv \sum_{d|n}d^{11} \ \bmod \ 691:$

\begin{verbatim}
gp > (Delta-h12)/691
%10 = -3*q^2 - 256*q^3 - 6075*q^4 - 70656*q^5 - 525300*q^6
 - 2861568*q^7 - 12437115*q^8 - 45414400*q^9
 - 144788634*q^10 - 412896000*q^11 - 1075797268*q^12
 - 2593575936*q^13 - 5863302600*q^14 - 12517805568*q^15
 - 25471460475*q^16 - 49597544448*q^17
 - 93053764671*q^18 - 168582124800*q^19 + O(q^20)
\end{verbatim}

\noindent
{\bf Применение модулярных форм к проблемам  теории чисел:}

\hskip-1cm
\begin{tabular}{cccc} 
\begin{tabular}{l}
Производящая\\ 
функция \\
$f=\sum_{n=0}^\infty a_n q^n$ \\
$\in {\mathbb C}[[ q]]$\\
для арифметической \\
функции $n\mapsto a_n$, \\
например $a_n=p(n)$
\end{tabular}
&
$\rightsquigarrow$
\begin{tabular}{l}
Выражение через\\
модулярную форму,\\
например \\
$\ds\sum_{n=0}^\infty p(n) q^n$\\
$=(\Delta/q)^{-1/24}$
\end{tabular}
&
$\rightsquigarrow$ 
\begin{tabular}{l}
Число\\
(ответ)
\end{tabular}
\\ 
\begin{tabular}{l}
{\bf Пример 1} (см. \rm \cite{Chand70}): \\ 
(Харди-Рамануджан)
\end{tabular} 
&
$\uparrow$
&
$\uparrow$
\\
\tiny
\begin{tabular}{l}
$\ds p(n)=\frac{e^{\pi \sqrt{2/3({n-1/24})}}} {4\sqrt{3}\lambda_n^2}$ \\  
$+O(e^{\pi \sqrt{2/3}\lambda_n}/ \lambda_n^3),$\\
$\lambda_n=\sqrt{n-1/24}$,
\end{tabular}
&
\begin{tabular}{l}
Хорошие базисы \\
конечномерность\\
много соотношений \\ 
и тождеств
\end{tabular}
&
\begin{tabular}{l}
Значения \\
$L$-функций, \\ 
сравнения, \\
\dots
\end{tabular}
\end{tabular}

\noindent
{\bf Пример 2}
(см. в \cite{Ma-Pa05}): теорема Ферма-Уайлза, гипотеза Бёрча-Суиннертона-Дайера, \dots

\section{\bf Семейства  модулярных форм Колмана}

{\bf Семейства  модулярных форм Колмана переменного веса $k\ge 2$}

\begin{tabular}{ll}
\begin{tabular}{l}
\begin{tabular}{l}
$\ds
k{}\mapsto f_{k{}}= \sum_{n=1}^\infty a_n(k{})q^n$\\
$\in \Qb[\![ q]\!] \subset {\mathbb C}_p[\![ q]\!]$\\ 
\end{tabular}
\\
\begin{tabular}{l}
Модельный пример\\
$p$-адического семейства:  
\end{tabular}
\\
\begin{tabular}{l}
{\it ряд Эйзенштейна веса} $k$\\
$\ds a_n(k{}) = \sum_{d|n}d^{k{}-1}, f_{k{}}=E_{k{}}$\\ 
$a_p(k)=1+p^{k-1}\Rightarrow$ Многочлен\\
Гекке (Эйлеровский множитель):\\
$1-a_pX+p^{k-1}X^2$\\ 
$=(1-X)(1-p^{k-1}X)$  
\end{tabular}
\end{tabular}
&
\hskip-0.5cm
\begin{tabular}{|l}
\begin{tabular}{l}
1) функции $k\mapsto a_n(k)$ \\
({\it коэффициенты ряда}) и\\ 
$p$-параметр Сатаке $k\mapsto \alpha_p^{(1)}(k{})$\\
являются $p$-адическими\\ 
аналитическими при $(n,p)=1$:\\
$1-a_pX+\psi(p)p^{k-1}X^2$\\ 
$=(1-\alpha_p^{(1)}(k)X)(1-\alpha_p^{(2)}(k)X)$  
\end{tabular}\\ 
\begin{tabular}{l}\\
2) наклон $\ord_p(\alpha(k))=\si > 0$  \\
{\it постоянен и положителен}
\end{tabular}
\end{tabular}
\end{tabular}

\subsection{\it Построение семейств Колмана}

{\bf Колман и Хида установили, что все классические модулярные формы живут в 
$p$-адических аналитических семействах}, см. \cite{Hi86}, \cite{CoPB}.
Комьютерная программа на PARI для вычисления таких семейств дана в \cite{CST98}
(см. также на Web-page of W.Stein, Modular Forms Database 
{\tt http://modular.fas.harvard.edu/})

{\bf Теория Банаховых модулей (Колман):}

\begin{tabular}{ll}
\begin{tabular}{l}
\noindent
$\bullet$
Оператор $U$ действует как \\
вполне непрерывный оператор \\ 
на Банахвом $\Ar$-подмодуле \\
${\cal M}^\dagger (N p^{v}; \Ar)$ \\ 
$\subset \Ar[\![ q]\!]$ (т.е. $U$ -- это предел \\ 
конечномерных  операторов)
\end{tabular} &
\begin{tabular}{l} 
$\Longrightarrow$ существует \\
{\it определитель Фредгольма} \\ 
$P_U(T)$\\ 
$=\det(Id-T\cdot U)\in \Ar[\![ T]\!]$\\
\\ 
\end{tabular}
\end{tabular}

\begin{tabular}{ll}
\begin{tabular}{l}
\noindent
$\bullet$
построен вариант \\
{\it теории Рисса}:\\
для любого обратного корня\\
$\al\in \Ar^*$ ряда $P_U(T)$ существует\\
собственная функция $g$, $Ug=\al g$\\
такая что все значения
\end{tabular} & 	
\begin{tabular}{l}
$ev_{k}(g)\in {{\mathbb C}_p}[\![ q]\!]$ \\ 
являются классическими \\
параболическими формами \\ 
для всех $k$ весов в некоторой \\ 
окрестности ${\cal B}\subset {X}$ \\
(см. в \cite{CoPB})
\end{tabular}
\end{tabular}

\section{\bf Тройные произведения Гарретта}

\subsection{\it Тройки (примитивных) модулярных форм}

$f_j(z)=\sum_{n=1}^\infty a_{n,j}q^n\in \Sr_{k_j}(N_j, \psi_j),\  (j=1, 2,3)$
весов $k_1, k_2, k_3$, кондукторов $N_1, N_2, N_3$, и характеров 
$\psi_j \bmod N_j$, $N:=  \lcm(N_1,N_2,N_3)$. 

Пусть $p$ простое число, $p\nmid N$.

Здесь 
$f_j\in \Qb[\![q]\!]\buildrel i_p\over\hookrightarrow{\mathbb C}_p[\![q]\!]$
при фиксированном вложении
$\Qb\buildrel i_p\over\hookrightarrow{\mathbb C}_p$, ${\mathbb C}_p={\widehat\Qb}_p$
поле Тэйта.
 
Рассматриваются только сбалансированные тройки весов:
\begin{eqnarray}\label{ptr0-1.3}
k_1 \ge k_2\ge k_3\ge 2,  \mbox{ и } k_1 \le k_2+ k_3-2 
\end{eqnarray}

\subsection{\it Тройное произведения Гарретта}

Это некоторое роизведение Эйлера степени 8:
\begin{align}\label{ptr0-1.1} &
 L(f_1\otimes f_2\otimes f_3, s, \chi) = 
\prod_{p\nmid N}
 L((f_1\otimes f_2\otimes f_3)_p, \chi(p) p^{-s}), 
\end{align}
\begin{align} \label{ptr1.2} &
\mbox{ где } 
 L((f_1\otimes f_2\otimes f_3)_p, X)^{-1}  = 
\\ 
& \det \left(1_8 -  X\mat \al^{(1)}_{ p,1}, 0, 0, \al^{(2)}_{ p,1}, \otimes
 \mat \al^{(1)}_{ p,2}, 0, 0, \al^{(2)}_{ p,2}, \otimes\mat 
\al^{(1)}_{ p,3}, 0, 0, 
\al^{(2)}_{ p,3},
 \right) \nonumber\\ 
&= \prod_{\eta}(1-\al_{ p,1}^{(\eta(1))}\al_{ p,2}^{(\eta(2))}
\al_{ p,3}^{(\eta(3))}X) \ \ \
\nonumber\\ 
&=(1-\al_{ p,1}^{(1)}\al_{ p,2}^{(1)}
\al_{ p,3}^{(1)}X)(1-\al_{ p,1}^{(1)}\al_{ p,2}^{(1)}
\al_{ p,3}^{(2)}X)\cds (1-\al_{ p,1}^{(2)}\al_{ p,2}^{(2)}
\al_{ p,3}^{(2)} X), 
\nonumber
\end{align}
произведение берется по всем 8 отображениям
$\eta: \{1, 2, 3 \}\to \{1, 2\}$.

\subsection{\it Критические значения и функциональное уравнение}

Нормализованная $Л$ функция (см. \cite{De79}, \cite{Co}, \cite{Co-PeRi}),
имеет вид:
\begin{align}\label{ptr6.4} & 
\La (f_1\otimes f_2\otimes f_3, s, \chi) = 
\\ &
\Ga_{\mathbb C}(s)
\Ga_{\mathbb C}(s-k_3+1)\Ga_{\mathbb C}(s-k_2+1)\Ga_{\mathbb C}(s-k_1+1)
 L(f_1\otimes f_2\otimes f_3, s, \chi),\nonumber
\end{align}

где $\Gamma _{{{\mathbb C}}}(s) = 2(2\pi )^{-s}\Gamma (s)$.
Гамма-множитель определяет критические значения $s=k_1, \cdots, k_2+ k_3-2$ 
для $\La(s)$.

{\bf  Функциональное уравнение} для $\La(s)$ имеет вид:
$$
s\mapsto k_1+k_2+k_3-2-s.
$$

\subsection{\it Метод: интегральное представлениые Гарретта}

{\bf Описание метода}

$\bullet$ 
вариант {\bf интегрального представления Гарретта} тройной $L$-функции в виде: для
$r=0, \cdots, k{}_2+k{}_3-k{}_1-2$, \\

$\ds \La (f_{1,k{}_1}\otimes f_{2,k{}_2}\otimes f_{3,k{}_3}, k{}_2+k{}_3-r, \chi)
 = $
\\

$\ds
\mathop{\int\int\int}_{\left(\Gamma_0(N^2p^{2v})\backslash {\HH}\right)^3}
\overline{{\tilde f_{1,k{}_1}}(z_1){\tilde f_{2,k{}_2}}(z_2){\tilde f_{3,k{}_3}}(z_3)}
{\cal E}({z}_1,{z}_2,{z}_3;-r, \chi)\prod_j ({dx_jdy_j\over y_j^2})
$

\noindent
где $\tilde f_{j,k_j}=:f_{j,k_j}^0$ это собственная функция сопряженного оператора Аткина
$U_p^*$ v ${\cal M}_{k_j}(Np, \psi_j)$, $f_{j,k_j, 0}$ -- собственная функция оператора
Аткина $U_p$, 
${\cal E}({z}_1,{z}_2,{z}_3;-r, \chi)\in {\cal M}_T(N^2p^{2v})
={\cal M}_{k_1{}}(N^2p^{2v}, \psi_1)\otimes {\cal M}_{k_2{}}(N^2p^{2v}, \psi_2)\otimes 
{\cal M}_{k_3{}}(N^2p^{2v}, \psi_3)$
тройная модулярная форма веса $(k{}_1,k{}_2, k{}_3)$, с фиксированным тройным характером
$(\psi_1, \psi_2, \psi_3)$.

\medskip
{\bf Тройная модулярная форма} ${\cal E}({z}_1,{z}_2,{z}_3;-r, \chi)$
строится из некоторого  почти голоморфного ряда {\it Зигеля-Эйзенштеина}
$F_{\chi, r}= G^{\star}({z}, -r; k{}, {(Np^v)}^2,\bpsi)$ рода три, веса 
$k{}=k{}_2+k{}_3-k{}_1$, и характера $\bpsi=\chi^2\psi_1\psi_2\overline\psi_3$. 
Для этого к ряду $F_{\chi, r}$ применяется:

{\it оператор скручивания Бёхерера} введённый в \cite{Boe-Sch},

{\it дифференциальный оператор Ибукиямы} (см. \cite{Ibu}, \cite{BSY}).

\medskip
{\bf Также используется:}

{\bf
$\bullet$
Теория $p$-адического интегрирования} 
со значениями в Банаховых ${\cal A}$-модулях ${\cal M}_T({\cal A})$ 
тройных  модулярных форм над $p$-адической Банаховой алгеброй ${\cal A}$. 
В этой теории строятся меры на группе $Y= (\Z/N\Z)^* \times\Z_p^*$ 
с использованием элементов ${\cal E}(-r, \chi)$ модуля ${\cal M}_T({\cal A})$. 

{\bf 
$\bullet$
Спектральная теория  тройного оператора  Аткина $U=U_ {p,T}$}
позволяет вычислить интеграл через проекцию $\pi_\lambda$ модуля 
${\cal M}_T({\cal A})$ на его $\lambda$-часть
${\cal M}_T({\cal A})^\lambda$.

Доказывается что $U$ -- это вполне непрерывный 
${\cal A}$-линейный оператор (т.е. предел конечномерных операторов),
и проекция $\pi_\lambda$ 
существует по общему результату  Серра и Колмана, см. \cite{CoPB}, \cite{SePB}.

\section{\bf Основной результат: $L$-функция четырёх переменных}

Обозначим $\chi$ переменный характер Дирихле ${}\bmod Np^v, v\ge 1$, 
и пусть $k_j{}$-переменные веса в {\bf пространстве $p$-адических весов 
$X=X_{Np^v}=\Hom_{cont}(Y, {\mathbb C}_p^*)$, $Y=(\Z/N\Z)^*\times\Z_p^*$}
($p$-адическое аналитическое пространство).
{\bf Для $r\in\Z$ точка $(r, \chi)\in X$ даётся гомоморфизмом 
$(y_1, y_2)\mapsto \chi(y_1)\chi(y_2\bmod p^v)y_2^r$.}

Символ $\ang{g,h},$ 
обозначает нормализованное скалярное произведение Петерсона модулярных форм.

{\bf Теорема.}
{\it
{\sf 1)} Функция
$\ds \Lr_{\underline f}: (s, k{}_1, k{}_2, k{}_3)\mapsto \frac {\ang{\underline f^0, 
{\cal E}(- r, \chi)},} { \ang{\underline f^0, \underline f_0},}$
зависит $p$-адически аналитически от четырех переменных
$(\chi\cdot y_p^r, k{}_1, k{}_2, k{}_3) \in X\times \Br_1\times\Br_2\times\Br_3$;
}

{\it
{\sf 2)} 
Для всех $p$-адических весов $({k{}_1, k{}_2, k{}_3})$ 
в некоторой  $p$-адической окрестности $\Br= \Br_1\times\Br_2\times\Br_3$,
с условием $k{}_1 \le k_2{} + k_3{} -2$, имеем: значения в точках 
$s=k_2{} + k_3{} -2-r$ совпадают с нормализованными критическими значениями
$$
L^*({f_{1,k_1{}}}\otimes {f_{2,k_2{}}}
\otimes {f_{3, k_3{}}}, k_2{}+k_3{}-2-r, \chi)  \ \ (r=0, \cdots, k_2{}+k_3{}-k{}_1-2), 
$$
для характеров Дирихле $\chi\bmod Np^v, v\ge 1$, с $Np$-полным кондуктором.}

{\it
{\sf 3)} 
Положим $H= [2{\rm ord}_p(\lambda)]+1$. Для всех весов 
$({k{}_1, k{}_2, k{}_3})\in\Br$ и $x=\chi\cdot y_p^r$
функция
$$
x\longmapsto 
\frac {\ang{\underline f^0, {\cal E}(- r, \chi)},}{ \ang{\underline f^0, \underline f_0},}
$$
продолжается до $p$-адической аналитической функции типа 
$o(\log^H(\cdot))$ переменной $x\in X$.}

\subsection{\it Набросок доказательства}

{\sf 1)} 
$\bullet$
Для любого веса $(k{}_1, k{}_2, k{}_3)$ используем равенство
${\ang{\underline f^0, {\cal E}(- r, \chi)},}={\ang{\underline f^0, 
\pi_\lambda({\cal E}(- r, \chi)}),}$
, которое выводится из
$\Im  \pi_\lambda=\Ker(U_T-\lambda I)^n, \Ker  \pi_\lambda=\Im(U_T-\lambda I)^n$.

Отсюда следует: ${\ang{\underline f^0, {\cal E}(- r, \chi)},}$ зависит
$p$-адически аналитически от $({k_1, k_2, k_3})\in\Br=\Br_1\times\Br_2\times\Br_3$, 
где ${\cal A}= {\cal  A}(\Br_1\times\Br_2\times\Br_3)$ -- $p$-адическая Банахова 
алгебра аналитических  функций на $\Br$.

$\bullet$
Для любого веса 
$(k{}_1, k{}_2, k{}_3)$ 
скалярное произведение 
$\ang{\underline f^0, {\cal E}(- r, \chi)},$
дается первой координатой вектора
$\pi_\lambda({\cal E}(- r, \chi))$ 
в любом ортогональном базисе модуля
${\cal M}^\lambda({\cal A})$
содержащем
$\underline f^0$ 
относительно алгебраического произведения Петерсона--Хиды 
$\ang{g, h},_a=\ang{g^\rho|\mat 0, -1, Np, 0, , h},$
(расширенного на модуль 
${\cal M}^\lambda({\cal A})$
по трилинейности).

{\bf Ответ:}

Выберем (локальный) базис $\ell^1, \cdots, \ell^n$ 
дающийся некоторыми тройными коэффициентами Фурье двойственного 
$\Ar$-модуля (локально свободного конечного ранга) 
${\cal M}^\lambda({\cal A})^*$. Используем обозначение
$\ell(h)=
\frac {\ang{\underline f^0, h},}{ \ang{\underline f^0, \underline f_0},}$. 

Геометрический смысл $\ell$: первая координата $\Ar$-базиса собственных функций 
операторов Гекке $T_q$ для всех $q\nmid Np$, с первым базисным элементом 
$f_0\in {\cal M}^\lambda({\cal A})$. 
Такой базис несложно построить с помощю действия $T_q$ на степенных рядах.

Кроме того, $\ell=\beta_1\ell^1+ \cdots+ \beta_n\ell^n$, и 
$\beta_i=\ell(\ell_i)$, причем $\ell_i$ обозначает двойственный базис модуля  
${\cal M}^\lambda({\cal A})$: $\ell^j(\ell_i)=\delta_{ij}$. Отсюда 
$\beta_i=\ell(\ell_i)=
\frac {\ang{\underline f^0, l^i},}{ \ang{\underline f^0, \underline f_0},}\in {\cal A}$.

Для любого веса $(k{}_1, k{}_2, k{}_3)$ имеем
$$
\ell({\cal E}(- r, \chi))=
\beta_1\ell^1({\cal E}(- r, \chi))+ \cdots+ \beta_n\ell^n({\cal E}(- r, \chi))
$$
где $\beta_i=\ell(\ell_i)\in \Ar$, причём  $\ell^i({\cal E}(- r, \chi))\in {\cal A}$ 
некоторые тройные коэффициенты Фурье функции ${\cal E}(- r, \chi)$.

Отсюда вытекает  1) поскольку все  коэффициенты Фурье
функции ${\cal E}(- r, \chi)$ лежат в $\Ar$.

{\bf Заключение}:
$\beta_i\ell_i({\cal E}(- r, \chi))$ лежат в $\Ar$ и
корректно определены для всех $p$-адических весов, не только
для классических весов $(k_1, k_2, k_3)$, откуда вытекает  1), 
что и даёт интерполяцию по переменным $(k_1, k_2, k_3)$.

Отсюда видно, что существует функция
$\tilde{\cal E}(- r, \chi)\in {\cal M}^\lambda({\cal A})$, 
такая, что 
$\ell({\cal E}(- r, \chi))=\ell(\pi_\lambda({\cal E}(- r, \chi))=\ell(\tilde{\cal E}(- r, \chi))$ 
(для любого веса $(k_1, k_2, k_3)$).

Для доказательства {\sf 2)},  {\sf 3)}, изучается зависимость от характера
$x=\chi\cdot y_p^r$ и теория $p$-адического интегрирования. 

Строятся $p$-адические меры $\mu$ со значениями в модуле 
${\cal M}^\lambda({\cal A})$.
Для этого мы доказываем сравнения для коэффициентов 
Фурье тройной модулярной формы  
${\cal E}({z}_1,{z}_2,{z}_3;-r, \chi)$, 
как в работах \cite{PaMMJ}, \cite{PaTV}.

Функция $\Lr$ строится как преобразование Меллина $\Lr(x)=\Lr_\mu(x)$ 
меры, которое всегда аналитически зависит от характера $x$.

Эти меры определены сначала лишь на ``пробных функциях'' $x=\chi\cdot y_p^r$, 
но они допускают каноническое продолжение на все локально-аналитические функции. 
\qed

\bibliographystyle{plain}

\end{document}